\documentclass[graybox]{svmult}

\usepackage{type1cm}        
\usepackage{makeidx}         
\usepackage{graphicx}        
\usepackage{multicol}        
\usepackage[bottom]{footmisc}

\usepackage{newtxtext}       %
\usepackage[varvw]{newtxmath}       

\usepackage{multirow}
\usepackage{diagbox}

\usepackage{booktabs}

\usepackage{float} 
\usepackage{bm}
\usepackage{subcaption}

\makeindex             
\def \a  {\alpha}

\def \e  {\varepsilon}

\def \l  {\lambda}
\def \om {\omega}
\def \Om {\Omega}

\def \P {\mathbb{P}}
\def \calT {\mathcal{T}}

\def \calA {\mathcal{A}}

\def \t  {\tau}

\def \del {\nabla}
\def \div {\nabla\cdot}

\def \p  {\partial}
\def \R  {\mathbb{R}}
\def \N  {\mathbb{N}}

\def \bmu {\bm{u}}
\def \bmv {\bm{v}}
\def \Id {\bm{I}}

\begin{document}

\title*{A History-dependent Dynamic Biot Model}
\author{Jakob S. Stokke\orcidID{0000-0002-2046-1475} and\\ Morten Jakobsen\orcidID{0000-0001-8861-1938} and\\ Kundan Kumar\orcidID{0000-0002-3784-4819} and\\ Florin A. Radu\orcidID{0000-0002-2577-5684}}
\authorrunning{Jakob S. Stokke et al.}
\institute{
Jakob S. Stokke \and Morten Jakobsen \and Kundan Kumar \and Florin A. Radu \at Center for Modeling of Coupled Subsurface Dynamics, University of Bergen, Bergen, Norway, 
\email{jakob.stokke@uib.no, morten.jakobsen@uib.no, kundan.kumar@uib.no,\\ florin.radu@uib.no}}

\maketitle

\abstract*{In this work, we consider a fully dynamic Biot model that includes memory effects due to evolving permeability. Time integrals are used to account for the change in structure. We propose a fully discrete scheme for this model, extending the fixed-stress splitting scheme. We use finite elements in space and a backward Euler discretization in time. The performance of the method is demonstrated through numerical experiments.}

\abstract{In this work, we consider a fully dynamic Biot model that includes memory effects due to evolving permeability. Time integrals are used to account for the change in structure. We propose an iterative splitting scheme for this model, extending the fixed-stress split for the quasi-static Biot. We use finite elements in space and a backward Euler discretization in time. The performance of the method is demonstrated through a numerical experiment. }

\section{Introduction}
\label{sec:intro}

There is a strong interest in coupling mechanics and flow in porous media due to many relevant societal applications. For example, there are several biomedical engineering, geophysical, and environmental applications, including heart perfusion, geothermal energy extraction, and $CO_2$ storage. The mathematical models for flow in deformable porous media consist of coupled (possibly nonlinear) partial differential equations. These models, based on Biot equations, are typically impossible to solve analytically and very challenging numerically. A lot of effort was made in the last decade for developing robust and efficient splitting schemes, but mostly for the quasi-static Biot model, see e.g. \cite{MikelicAndro2013Coic, both2017etal, Borregales2018,StorvikErlend2019Otoo}.


However, when the ratio between the characteristic time and the characteristic domain time scale is large, a fully dynamic Biot model must be considered. In this work we consider the Biot-Allard model as derived by homogenization in  \cite{MikelicBiotAllard}. Specific for this model is that it includes an acceleration term in the mechanics equation and memory effects expressed as convolution integrals. The system without memory effects, which consists of two linear fully coupled equations, one hyperbolic and one parabolic,  was already studied in \cite{ShowalterR.E.2000DiPM, BauseMarkusSpaceTime}. A splitting scheme based on the undrained-split \cite{KimJihoon2011SAaE} was derived and analyzed in \cite{BauseMarkus2020ICfF} following a gradient flow approach \cite{both2019gradient}. 




In this paper, we propose a fixed-stress \cite{KimJihoon2011SAaE} type splitting scheme for the fully dynamic Biot-Allard model including memory effects (it includes a history-dependent permeability). For the discretization in time we use an implicit Euler scheme, and in space we use Galerkin finite elements. A numerical experiment is performed to study the performance of the scheme.

The paper is structured as follows. In the next section, we present the mathematical model. Section \ref{sec_numscheme} concerns the discretization and the splitting scheme. In Section \ref{sec_examples} we have two numerical examples with a manufactured solution. The paper ends with a short concluding section.

\section{A History-dependent Dynamic Biot Model}\label{sec_model}

We consider flow in a saturated porous medium $\Om\subset\R^{d}$, with $d\in\{2,3\}$ being the spatial dimension, which is elastic, homogeneous, and isotropic, and has a Lipschitz continuous boundary $\partial \Om$. 
The history-dependent dynamic Biot model we consider is a simplified Biot-Allard model: Find $(\bmu,p)$ in $\Om\times [0,T]$ such that
\begin{subequations}\label{papereqs:fullydynamicBiot}
		\begin{align}
			\frac1M\p_t p+\alpha\p_{t}\div\bmu+\div\left[\int_{0}^{t}\calA\left(\frac{(t-\zeta)\eta}{l^2\rho_f}\right)\left(-\frac{1}{\rho_f}\del p(x,\zeta)\right)\right]d\zeta&=f_2,\label{eq:biotallardflow}\\
            \rho\p_{t}^2\bmu-\div\left[2\mu\e(\bm{u})+\lambda\div\bmu\Id-\alpha p\Id\right]\quad\quad&\quad\quad\nonumber\\+\p_t\int_{0}^{t}\calA\left(\frac{(t-\zeta)\eta}{l^2\rho_f}\right)\left(-\del p(x,\zeta)\right)d\zeta&=\textbf{f}_1,\label{eq:biotallardmechanincs}
		\end{align}
  with initial and boundary conditions
\begin{align}
      \bmu(0)=\bmu_0,\quad\p_t\bmu(0)=\bmu_1, \quad p(0)=p_0,\quad &\mbox{ in }\Om,\\
             \bmu=0, \quad p=0,\quad &\mbox{ on }\p\Om\times(0,T].
\end{align}
	\end{subequations}
 Here, $\bmu$ denotes the solid phase displacement and $p$ the pressure. Furthermore, $\alpha$ is the Biot coefficient, $\rho=\rho_f\varphi+\rho_s(1-\varphi)$ is the mass density with $\varphi$, $\rho_s$ and $\rho_f$ being the porosity, solid and fluid density respectively, $\lambda$ and $\mu$ are the Lamé parameters, $\frac1M$ is the compressibility constant, $\eta$ is the pore fluid viscosity, $l$ the typical pore size, $\bm{I}\in\R^{d\times d}$ is the identity matrix and $\e(\bmu)=\frac12(\del\bmu+\del\bmu^{T})$ the linearized strain tensor. The right-hand side functions of \eqref{papereqs:fullydynamicBiot}, $\bm{f}_1$ and $f_2$, are an external force density and a volumetric source, respectively. The changes in permeability are included in $\calA(\cdot)$ the dynamic permeability tensor. It can be computed by solving a Stokes system at the pore-scale (see \cite[p.15]{MikelicBiotAllard}). In this work we will assume that the solution to the Stokes problem is known. We refer to the recent paper \cite{stokke2024a} for the well-posedness of the dynamic Biot model with memory effects.

\section{A Fully Discrete Numerical Scheme}\label{sec_numscheme}
We use functional analysis notations throughout our paper. The space of Lebesque measurable and square-integrable functions is denoted by $L^{2}(\Om)$. By $\left\langle\cdot,\cdot\right\rangle$ and $\|\cdot\|$ we denote the $L^2(\Om)$ scalar product and norm. Further $H^{1}(\Om)$ is the space with $L^{2}$-functions having weak first order derivatives in $L^{2}(\Om)$ and $H^{1}_0(\Om):=\left\{f\in H^{1}(\Om)\,\, |f=0 \mbox{ on }\p\Om\right\}$.

\subsection{Discretization of Convolution Integrals}

To discretize \eqref{papereqs:fullydynamicBiot}, we first need to deal with the convolution terms. We neglect the scaling factor  $\eta/(l^{2}\rho_f)$ for the ease of the presentation. We also assume that a part of the history of the porous medium is known before the initial time, i.e., $\calA(-\t)$ for a $\tau>0$. The convolutions are approximated using the trapezoid rule
\begin{align*}
    &\int_{0}^{t_n}\calA\left(t_n-\zeta\right)\left(-\del p(x,\zeta)\right) d\zeta
    \approx -\frac12\big(\calA\left(0\right)\del p(x,t_n)\\
    &+2\sum_{i=1}^{n-1}\calA\left(t_n-\zeta_i\right)\del p(x,\zeta_i)+\calA\left(t_n\right)\del p(x,0)\big)\t := -\sum_{i=0}^{n}\t \omega_i \calA^{i}\del p^{i},
\end{align*}
where $\t=\zeta_i-\zeta_{i-1}$ and $\omega_i$ is the corresponding quadrature weight.
In addition, we need to approximate the time derivative of the convolution
\begin{align*}
    \left(\p_t\int_{0}^{t}\calA\left(t-\zeta\right)\del p(x,\zeta)d\zeta\right)\Big|_{t=t_{n}}=&\, \calA(0)\del p(x,t_n)\\
    &+\int_0^{t_n}\p_t \calA\left(t_n-\zeta\right)\del p(x,\zeta)d\zeta,
\end{align*}
which, by use of the backward difference approximation and the trapezoid rule, can be approximated using

\begin{align*}
   \approx&\,\calA(0)\del p(x,t_n)+\sum_{i=0}^{n}\om_i\left[\calA(t_n-\zeta_i)-\calA(t_n-\t-\zeta_i)\right]\del p(x,\zeta_i) \\ :=&\,\calA^{0}\del p^n +\sum_{i=0}^{n}\omega_i \Delta\calA^{i}\del p^{i},\nonumber 
\end{align*}
where $\Delta \calA^{i} = \calA^{i}-\calA^{i-1}$.

\subsection{Discretization in Space and Time}
Let $\calT_h$ be a regular decomposition of the domain $\Om$,  with $h$ denoting the mesh diameter. For the spatial discretization, we use the Taylor-Hood elements for the displacement and pressure
\begin{align*}
    V_h=&\left\{\bmv_h\in \left[H_0^{1}(\Om)\right]^d\, \big|\,\forall\, K\in\calT_h,\, \bmv_h|_K\in [\P_2]^d\right\},\\
    Q_h=&\left\{q_h\in H_0^{1}(\Om)\, \big|\,\forall\, K\in\calT_h,\, q_h|_K\in \P_1\right\},
\end{align*}
where $\P_1$ and $\P_2$ denote the spaces of linear and quadratic polynomials, respectively.  We consider a uniform partition of the time interval $[0,T]$, i.e. $0<t_1<\cdots <t_N=T$ with constant time step size $\t=t_n-t_{n-1}$. The temporal derivatives in \eqref{papereqs:fullydynamicBiot} are approximated using finite differences,
\begin{align*}
    \p_t p(t_n)\approx \frac{p^n-p^{n-1}}{\t},\quad \p_t^{2}\bmu(t_n)\approx\frac{\bmu^{n}-2\bmu^{n-1}+\bmu^{n-2}}{\t^{2}}.
\end{align*}

We assume that the initial data is known and that the first time step has already been determined, then we can approximate the history-dependent dynamic Biot equations by the implicit Euler scheme:  For $n\geq 2$ given $\bmu^{n-1}_{h},\bmu^{n-2}_{h}$ and $p^{n-1}_{h}$ find $(\bmu_h^n,p_h^n)\in V_h\times Q_h$ which for all $(\bmv_h,q_h)\in V_h\times Q_h$ satisfies


\begin{subequations}\label{eqs:2}

\begin{align}
        \frac{1}{M}\left\langle p^{n}_h-p^{n-1}_h,q_h\right\rangle
        +  \left\langle\alpha\div\left(\bmu^{n}_h-\bmu^{n-1}_h\right),q_h\right\rangle\nonumber\\
         +\t^{2}\left\langle \frac{1}{\rho_f}\sum_{i=0}^{n} \omega_i \calA^{i}\del p^{i}_{h},\del q_h\right\rangle&=\t\left\langle f_2,q_h\right\rangle,\\
        \rho\left\langle \bmu^{n}_{h}-2\bmu^{n-1}_{h}+\bmu^{n-2}_{h}, \bmv_h\right\rangle + \t^{2}\lambda\left\langle \div\bmu_h^{n}, \div\bmv_h\right\rangle\nonumber\\
        +2\t^{2}\mu\left\langle \e(\bmu_h^{n}), \e(\bmv_h)\right\rangle-\t^{2}\a\left\langle p_h^{n},\div \bmv_h\right\rangle\nonumber\\-\t^{2}\left\langle\calA^0\del p^{n}_{h},\bmv_h\right\rangle- \t^{2}\left\langle \sum_{i=0}^{n}\omega_i \Delta\calA^{i}\del p^{i}_{h},\bmv_h\right\rangle&=\t^{2}\left\langle \bm{f}_1,\bmv_h\right\rangle.
    \end{align}
\end{subequations}
\subsection{Fixed-Stress Splitting Scheme}
\begin{subequations}\label{eqs:fixedstress}
Let $k\in\N$ denote the iteration index and note that in the history terms of \eqref{eqs:2}, we sum over $p^{i}_h$ including $p^{n}_h$ which we take to be the $k^{\rm th}$ iterate $p^{n,k}_{h}$. Now based on the fixed-stress scheme for the quasi-static Biot \cite{KimJihoon2011SAaE}, we can write a splitting scheme for the model \eqref{papereqs:fullydynamicBiot}. At any $k\geq 1$ and for a free-to-be-chosen stabilization parameter $L>0$, (we fix the value later, see remark \ref{remark1}) the scheme reads:

Given $\bmu^{n,k-1}_{h}\in V_h$ find $p^{n,k}_{h}\in Q_h$ such that
\begin{align}
        \frac{1}{M}\left\langle p^{n,k}_h-p^{n-1}_h,q_h\right\rangle+L\left\langle p^{n,k}_h-p^{n,k-1}_h,q_h\right\rangle&\nonumber\\
        +  \left\langle\alpha\div\left(\bmu^{n,k-1}_h-\bmu^{n-1}_h\right),q_h\right\rangle
         +\t^{2}\left\langle \frac{1}{\rho_f}\sum_{i=0}^{n} \omega_i \calA^{i}\del p^{i}_{h},\del q_h\right\rangle&=\t\left\langle f_2,q_h\right\rangle.
    \end{align}
    The second step is solving the mechanical part: Given $p^{n,k}_{h}\in Q_h$ find $\bmu^{n,k}_{h}\in Q_h$ such that
    \begin{align}
        \rho\left\langle \bmu^{n,k}_{h}-2\bmu^{n-1}_{h}+\bmu^{n-2}_{h}, \bmv_h\right\rangle + \t^{2}\lambda\left\langle \div\bmu_h^{n,k}, \div\bmv_h\right\rangle&\nonumber\\
        +2\t^{2}\mu\left\langle \e(\bmu_h^{n,k}), \e(\bmv_h)\right\rangle-\t^{2}\a\left\langle p_h^{n,k},\div \bmv_h\right\rangle&\nonumber\\ -\t^{2}\left\langle\calA^0\del p^{n,k}_{h},\bmv_h\right\rangle- \t^{2}\left\langle \sum_{i=0}^{n}\omega_i \Delta\calA^{i}\del p^{i}_{h},\bmv_h\right\rangle&=\t^{2}\left\langle \bm{f}_1,\bmv_h\right\rangle.
    \end{align}
\end{subequations}


\begin{remark}\label{remark1}
    The convergence of the scheme depends on the choice of the parameter $L$. A rigorous analysis was performed for the case of the linear quasi-static Biot in e.g. \cite{MikelicAndro2013Coic,both2017etal} or for nonlinear quasi-static Biot in e.g. \cite{Borregales2018, Kraus2023}. Finding an optimal $L$, in the sense of minimizing the number of iterations, is not an easy task \cite{StorvikErlend2019Otoo}, in practice $L=\frac{\alpha}{K_{dr}}$ works well, where $K_{dr}=\frac{2\mu}{d}+\l$.
\end{remark}

 \section{Numerical results}\label{sec_examples}
In this section, we will perform two numerical experiments to study the convergence of the proposed splitting scheme.  As a stopping criterion, we have used a similar relative $L_2$ norm as in \cite{StorvikErlend2019Otoo}, i.e. $\|p^{n, i}_{h}-p^{n,i-1}_{h}\|\leq \e_r\|p^{n, i}_{h}\|$, with $\e_r=10^{-9}$. 

The stabilization parameter $L$ we have used here corresponds to the one typically used in the fixed-stress scheme for the quasi-static Biot model, i.e. $L_{C}=\frac{\alpha^{2}}{K_{dr}}$, cf., e.g., \cite{KimJihoon2011SAaE}. The linear systems are solved using a direct solver. All numerical experiments were implemented using the open-source software DUNE \cite{DunePaper}.  In both examples we solve the Biot problem \eqref{papereqs:fullydynamicBiot} in the unit-square $\Om=(0,1)^{2}$ and with the paremeters $k_0=\alpha=M=\rho=\rho_f=\eta=1$, and $\mu=\l=10$.
\subsection{Example 1}

\begin{figure}[H]
\centering
 \begin{subfigure}[b]{0.48\textwidth}
   \includegraphics[width=\textwidth]{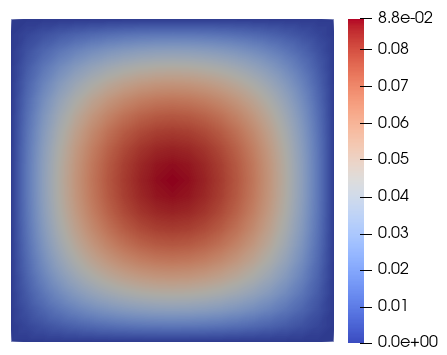}
   \caption{}
    \label{fig:subfig1}
  \end{subfigure}
  \hfill
\begin{subfigure}[b]{0.48\textwidth}
 \includegraphics[width=\textwidth]{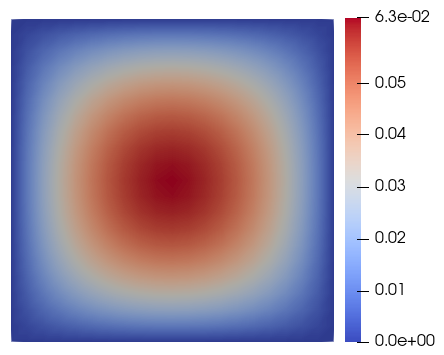}
 \caption{}
  \label{fig:subfig2}
\end{subfigure}
 \caption{Numerical solutions at $T=1$; (a) Displacement $|u_h|$, (b) Pressure $p$.}
\label{fig:ex1}
\end{figure}

 For the first example, the problem is solved on the time interval $I=[0,1]$. In all cases, the time step size is $\t=0.1$.
We let the dynamic permeability be
 \begin{equation}\label{eq:dynamicperm}
     \calA(t)=k_0+0.02k_0\sin(\pi t),
 \end{equation}
and construct source terms corresponding to 
 \begin{align*}
     u_1(x,y,t)=u_2(x,y,t)=t^{2}xy(1-x)(1-y),
 \end{align*}
 and 
 \begin{align*}
     p(x,y,t)=txy(1-x)(1-y).
 \end{align*}
On the boundary, we impose homogeneous Dirichlet boundary conditions for both the pressure and the displacement. The pressure and displacement magnitude at final time $T=1$ are displayed in Figure \ref{fig:ex1}.

\begin{figure}[H]
 \centering
 \begin{subfigure}[b]{0.48\textwidth}
   \includegraphics[width=\textwidth]{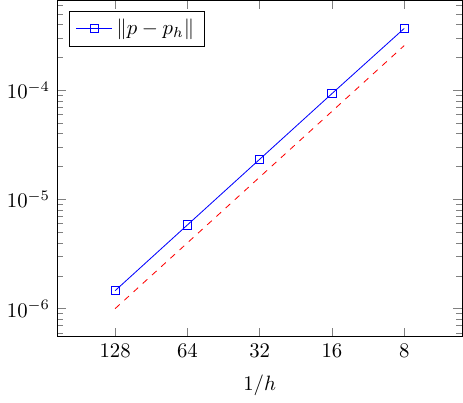}
 \end{subfigure}
 \hfill
 \begin{subfigure}[b]{0.51\textwidth}
   \includegraphics[width=\textwidth]{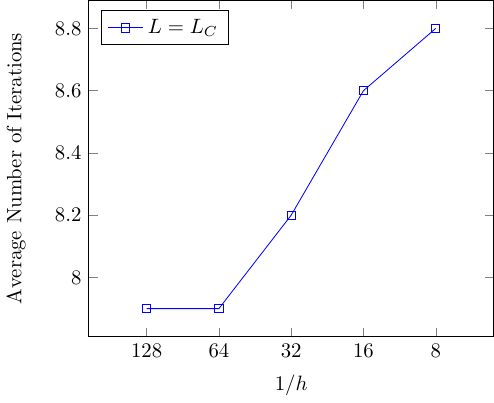}
\end{subfigure}
\caption{$L^{2}$ error of the pressure with a line (red) representing second-order convergence and the average number of iterations for different mesh sizes.}
\label{fig:convergence}
\end{figure}

Figure \ref{fig:convergence} shows the $L^{2}$ error of the pressure and the average number of iterations per time step for different mesh sizes. As the mesh is refined, the pressure converges with an order of 2. The average number of iterations is nearly mesh-independent, consistent with the theoretically established results for the quasi-static Biot model.  

\begin{figure}
     \centering
     \includegraphics[width=0.6\textwidth]{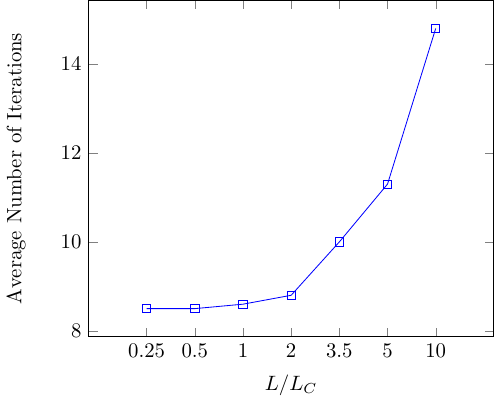}
     \caption{Average number of iterations per time step for different stabilization parameters with fixed mesh size $h=1/16$.}
     \label{fig:Lchange}
 \end{figure}

The average number of iterations per time step for different stabilization parameters with a fixed mesh size is displayed in Figure \ref{fig:Lchange}. In accordance with the quasi-static scheme, the stabilization parameter heavily influences the convergence rate. The fastest convergence is obtained for a smaller $L$ than the classical parameter.

 \subsection{Example 2}
\begin{figure}[H]
\centering
 \begin{subfigure}[b]{0.48\textwidth}
   \includegraphics[width=\textwidth]{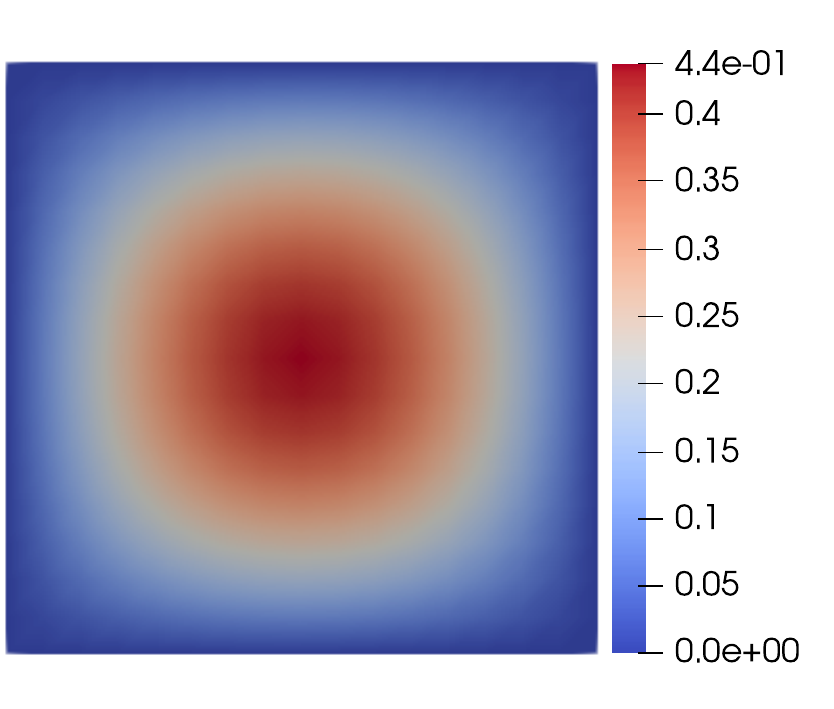}
   \caption{}
    \label{fig:subfig3}
  \end{subfigure}
  \hfill
\begin{subfigure}[b]{0.48\textwidth}
 \includegraphics[width=\textwidth]{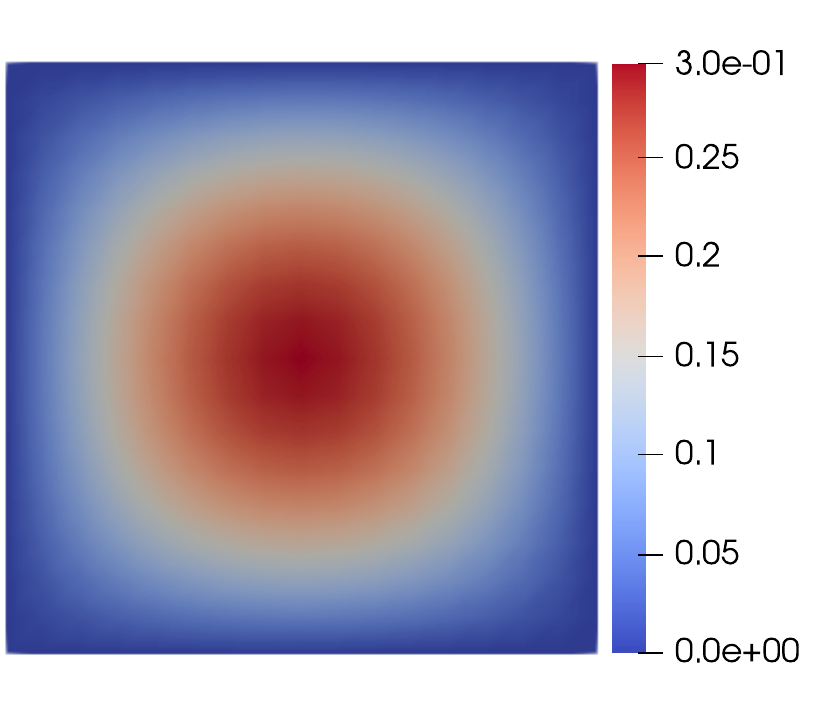}
 \caption{}
  \label{fig:subfig4}
\end{subfigure}
 \caption{Numerical solutions at $T=0.1$; (a) Displacement $|u_h|$, (b) Pressure $p$.}
\label{fig:ex2}
\end{figure}

For the second example, we consider the analytical solution
\begin{align*}
     u_1(x,y,t)=u_2(x,y,t)=p(x,y,t)=sin(\pi t)sin(\pi x)\sin(\pi y).
 \end{align*}
 The parameters are the same as in Example 1, and the dynamic permeability is given in  \eqref{eq:dynamicperm}.  We impose homogeneous Dirichlet boundary conditions for both pressure and displacement. Here, we solve the problem on the time interval $I=[0,0.1]$ with a variable time step size $\t$ and mesh size $h$. We have used the stabilization parameter $L=L_{c}$ in all computed solutions. The numerical solution at the final time $T=0.1$ is displayed in Figure \ref{fig:ex2}.

 \begin{figure}[H]
 \centering
 \begin{subfigure}[b]{0.48\textwidth}
   \includegraphics[width=\textwidth]{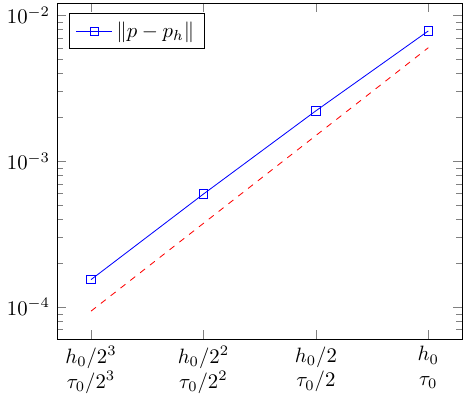}
 \end{subfigure}
 \hfill
 \begin{subfigure}[b]{0.51\textwidth}
   \includegraphics[width=\textwidth]{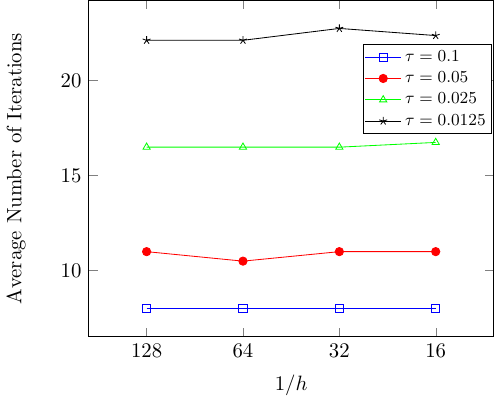}
\end{subfigure}
\caption{$L^{2}$ error of the pressure with a line (red) representing second-order convergence and the average number of iterations for different mesh sizes and time step sizes. Here $h_0=1/16$ and $\tau_{0}=0.1$.}
\label{fig:error}
\end{figure}
The $L^{2}$ error of the pressure as both the time step and mesh size are refined is shown in Figure \ref{fig:error}. We observe a convergence order of 2 for the pressure. Figure \ref{fig:error} also displays the average number of iterations with respect to mesh size for different time step sizes. Again, the number of iterations is nearly mesh-independent. However, it increases as the time step decreases. This is consistent with the contraction rate in theoretical estimates for the fixed-stress in the quasi-static case.
\section{Conclusions}

We considered a fully dynamic Biot model for flow and deformation in porous media. The model includes additional memory effects, containing a history-dependent permeability. We propose a fixed-stress type splitting scheme for efficiently solving this model. We performed two numerical experiments to study the convergence of the scheme.  It is not clear under which situations the memory effects become important. This will be explored under more realistic parameter choices in further work.

\begin{acknowledgement}
The authors acknowledge the support of the VISTA program, The Norwegian Academy of Science and Letters, and Equinor. The authors would like
to thank the reviewers for helpful comments to improve this work.
\end{acknowledgement}


\end{document}